\documentclass[11pt,a4paper]{article}
\usepackage{graphicx} 
\usepackage{amssymb}
\usepackage{amsmath}
\usepackage{latexsym, amsfonts}
\usepackage{xcolor}
\usepackage{dsfont}
\allowdisplaybreaks
\newcommand{\be}{\begin{equation}}
\newcommand{\ee}{\end{equation}}
\newcommand{\bea}{\begin{eqnarray}}
\newcommand{\eea}{\end{eqnarray}}
\newcommand{\bean}{\begin{eqnarray*}}
\newcommand{\eean}{\end{eqnarray*}}
\newcommand{\brray}{\begin{array}}
\newcommand{\erray}{\end{array}}
\newcommand{\ben}{\begin{equation}{nonumber}}
\newcommand{\een}{\end{equation}{nonumber}}

\newtheorem{dfn}{Definition}[section]
\newtheorem{thm}[dfn]{Theorem}
\newtheorem{lmma}[dfn]{Lemma}
\newtheorem{ppsn}[dfn]{Proposition}
\newtheorem{crlre}[dfn]{Corollary}
\newtheorem{xmpl}[dfn]{Example}
\newtheorem{rmrk}[dfn]{Remark}
\newcommand{\bdfn}{\begin{dfn}}
\newcommand{\bthm}{\begin{thm}}
\newcommand{\blmma}{\begin{lmma}}
\newcommand{\bppsn}{\begin{ppsn}}
\newcommand{\bcrlre}{\begin{crlre}}
\newcommand{\bxmpl}{\begin{xmpl}}
\newcommand{\brmrk}{\begin{rmrk}}
\newcommand{\edfn}{\end{dfn}}
\newcommand{\ethm}{\end{thm}}
\newcommand{\elmma}{\end{lmma}}
\newcommand{\eppsn}{\end{ppsn}}
\newcommand{\ecrlre}{\end{crlre}}
\newcommand{\exmpl}{\end{xmpl}}
\newcommand{\ermrk}{\end{rmrk}}

\newcommand{\IC}{\mathbb{C}}

\newcommand{\al}{\alpha}

\newcommand{\eps}{\epsilon}

\newcommand{\cla}{{\cal A}}
\newcommand{\clb}{{\cal B}}
\newcommand{\clc}{{\cal C}}
\newcommand{\cld}{{\cal D}}

\newcommand{\clg}{{\cal G}}
\newcommand{\clh}{{\cal H}}

\newcommand{\clk}{{\cal K}}
\newcommand{\cll}{{\cal L}}
\newcommand{\clm}{{\cal M}}
\newcommand{\cln}{{\cal N}}

\newcommand{\clq}{{\cal Q}}

\newcommand{\cls}{{\cal S}}

\newcommand{\clu}{{\cal U}}
\newcommand{\clv}{{\cal V}}
\newcommand{\clw}{{\cal W}}

\def\a*{{\cal A}_{h,*}}
\def\B{{\cal B}(h)}
\def\B1{{\cal B}_1(h)}
\def\b{{\cal B}^{\rm s.a.}(h)}
\def\b1{{\cal B}^{\rm s.a.}_1(h)}

\newcommand{\ot}{\otimes}

\newcommand{\raro}{\rightarrow}

\def \qed {$\Box$}

\begin{document}
\[
\]
\begin{center}
{\large {\bf  Linear coactions of discrete quantum groups on the circle }}\\
by\\
{\large Debashish Goswami{\footnote{ Partially supported by JC Bose National Fellowship given by SERB, Govt. of India.}}
},
{\large Suchetana Samadder{\footnote{ Support by a research fellowship from Indian Statistical Institute is gratefully acknowledged.}}}\\
{ Stat-Math Unit, Kolkata,}\\
{ Indian Statistical Institute,}\\
{ 203, B. T. Road, Kolkata 700 108, India}\\
{e mails: debashish\_goswami@yahoo.co.in, suchetana.smdr@gmail.com }

\end{center}

\begin{abstract}
For a (unital) $C^*$-algebra $\cla$, we construct a $C^*$-algebraic discrete quantum group (DQG) $\clq_{\rm aut}(\cla)$, coacting on $\cla$, which is a quantum generalization of ${\text Aut}(\cla)$ in the framework of discrete quantum groups, in the sense that any other coaction of a DQG on $\cla$ factors through the above coaction of $\clq_{\rm aut}(\cla)$.  We prove by an explicit calculation that if any Kac-type $C^*$-algebraic discrete quantum group $\mathcal{Q}$ has a `weakly faithful' coaction  on $C(S^1)$ which is `linear' in the sense that it leaves the space spanned by $\{ Z, \overline{Z} \}$ invariant, then  $\mathcal{Q}$ must be classical, i.e. isomorphic with $C_0(\Gamma)$ for some discrete group $\Gamma$. 
This parallels the well-known result of  non-existence of genuine compact quantum group symmetry obtained by the first author and his collaborators (\cite{sir_book} and the references therein).
\end{abstract}
\section{Introduction}
Quantum groups are well-known symmetry objects in mathematics and physics, originating from the pioneering work of Drinfeld and Jimbo (\cite{drinfeld}, \cite{jimbo}), in the algebraic setting and then also in the functional analytic framework by Woronowicz (compact quantum groups, \cite{woro1},\cite{woro2}), Van Daele (discrete quantum groups, \cite{dqg}), Vaes-Kustermans (locally compact quantum groups, \cite{vnqgp}), Baaj-Skandalis (multiplicative unitary , \cite{baaj_skandalis}), just to mention only a small fraction of the vast and expanding literature. It is natural to consider generalization of full automorphism or symmetry of some given mathematical structure (sets, algebras and so on) given by quantum groups of suitable kind. A lot of work has already been done by many authors in this direction too, beginning with Manin (\cite{manin_book}), Wang, Banica, Bichon (\cite{wang1}, \cite{wang2}, \cite{ban},\cite{bichon},\cite{bichon_banica}), just to name a few. The first author and his collaborators studied in details symmetry given by compact quantum groups on classical spaces, particularly smooth manifolds and proved (\cite{goswami_joardar}, \cite{goswami}) a non-existence result for genuine compact quantum group symmetry on compact connected smooth manifolds under the mild assumption of smoothness of the coaction giving the quantum  symmetry. However, it is only recently that quantum symmetry given by discrete and even locally compact quantum groups has gained attention among experts (\cite{voigt}, \cite{rollier_vaes}). It is natural to ask whether one has a similar non-existence result for discrete quantum group symmetry on compact connected smooth manifolds as well, or whether there is a chance of getting genuine quantum symmetry when we turn to discrete quantum groups. \\In order to investigate this question, one needs to have a formulation of a quantum version of ${\rm Aut}(\cla)$ in the discrete quantum group framework. We have been able to come up with such a notion of discrete quantum group of automorphisms of a unital $C^*$-algebra $\cla$, namely, $\clq_{\rm aut}(\cla)$, which has been proved in Theorem \ref{qaut_univ} generalizing a similar result proved by Voigt \cite{voigt} for the $C^*$-algebra $C_0(X)$ for an infinite countable set $X$.\\We then explore the existence of genuine discrete quantum group symmetry for a very elementary example, that of a circle, in this paper.
By explicit computations we prove that there is no genuine (noncommutative as $C^*$-algebra) discrete quantum group having a weakly faithful coaction on $C(S^1)$ which leaves the linear span of the canonical coordinate functions $Z, \overline{Z}$ invariant.
\section{Preliminaries: A quick review of compact and discrete quantum groups and their coactions on $C^*$-algebras}
A $C^*$-algebraic compact quantum group (CQG) is given by a unital $C^*$-algebra $\hat{\clq}$ with a unital $\ast$-homomorphism $\Delta$ from $\hat{\clq}$ to the minimal tensor product $\hat{\clq} \ot_{\min}  \hat{\clq}$ which is co-associative ($(\Delta \ot {\rm id}) \circ \Delta=({\rm id } \ot \Delta) \circ \Delta)$ and each of the linear spans ${\rm Span}(\Delta(\hat{\clq}) (1 \ot \hat{\clq}) )$ and ${\rm Span}(\Delta(\hat{\clq}) (\hat{\clq} \ot 1) )$ is norm dense in $\hat{\clq} \ot_{\rm min} \hat{\clq}$. A $C^*$-algebraic 
 discrete quantum group (DQG) is a non-unital (unless finite dimensional) $C^*$-algebra $\clq$ given by a direct sum of finite dimensional matrix algebras of the form $M_{d_i}, i \in I$ indexed by some set $I$, equipped with a coproduct  from $\clq$ into $\clm(\clq \ot_{\rm min} \clq)$ satisfying certain technical conditions as in \cite{dqg}, where for any $C^*$-algebra $\clc$, $\clm(\clc)$ denotes the 
  multiplier algebra of $\clc$. On a $C^*$-algebraic DQG there exists a bounded counit $\epsilon$ given by a central projection $P_{\mathds{1}}$ onto the one-dimensional direct summand coming from the unit corepresentation $\mathds{1}$. We also denote by $\epsilon$, the corresponding normal extension to the von Neumann algebra $\clq^{\rm von}=\prod_{i\in I}M_{d_i}$. In particular, given any von Neumann algebra $\clc^{\rm von}$, we have $({\rm id}\otimes \epsilon)$ as a well defined bounded normal $\ast$-homomorphism on $\clc^{\rm von}\overline{\otimes}\clq^{\rm von}$, where $\overline{\otimes}$ denotes the von Neumann algebraic tensor product.
  There is a natural duality between CQG and DQG inside the bigger category of locally compact quantum groups. More details can be found in the references already cited in the introduction. 
  
  In the von Neumann algebraic framework, a von Neumann algebraic CQG is a von Neumann algebra $\hat{\clq}^{\rm von}$ equipped with a normal unital $\ast$-homomorphism $\Delta$ from $\hat{\clq}^{\rm von}$ to $\hat{\clq}^{\rm von} \overline{\ot} \hat{\clq}^{\rm von}$
  and a faithful state $\tau$ which is bi-invariant in the sense that $(\tau \ot {\rm id}) \circ \Delta(a)=\tau(a)1=({\rm id } \ot \tau) \circ \Delta(a)$ for all $a \in \hat{\clq}^{\rm von}$. Indeed, $C^*$ and von Neumann algebraic definitions of CQG are essentially equivalent as the so-called reduced von Neumann algebra for  any reduced $C^*$-algebraic CQG with respect to the Haar state becomes a von Neumann algebraic CQG and any von Neumann algebraic CQG arises in this way from a unique reduced $C^*$-algebraic CQG. Similarly, we have a von Neumann algebraic DQG which is a von Neumann algebra $\clq^{\rm von}$ generated by finite dimensional matrix algebras $M_{d_i}, i \in I$, say along with a coassociative coproduct $\Delta$ from $\clq^{\rm von}$ to $\clq^{\rm von} \overline{\ot} \clq^{\rm von}$, which is normal, unital and also a pair of 
   faithful semifinite weight $\tau_l$ and $\tau_r$ which are left and right invariant respectively w.r.t. $\Delta$ in the sense of \cite{vnqgp}. We remark that there is an exact analogue of corepresentation theory and Peter-Weyl theory for CQG, as discussed in \cite{woro1}, \cite{woro2}.
   
    We call a $C^*$-algebraic or von Neumann algebraic compact quantum group to be of Kac type if the Haar state $\tau$ is tracial. In the dual picture a $C^*$-algebraic or von Neumann algebraic discrete quantum group is said to be of Kac type if the left and right Haar weights are equal. From the duality, one can easily infer that in both the $C^*$-algebraic and von Neumann algebraic settings a CQG is of Kac type if and only if the dual DQG is of Kac type.
   
   One also has analogue of action by a group in the realm of quantum groups. These are precisely given by co-actions on a $C^*$-algebra. We refer to \cite{neshveyev_tuset} for a definition of $C^*$-coaction of a $C^*$-algebraic locally compact quantum group $\clq$ on a unital $C^*$-algebra $\cla$, which is a coassociative, injective and non-degenerate $C^*$-algebra homomorphism $\delta:\cla\to \clm(\cla\otimes _{\rm min}\clq)$ such that the span density condition holds, i.e. $\delta(\cla)(1_{\cla}\otimes_{\rm min}\clq)$ is norm dense in $\cla\otimes_{\rm min}\clq$. However, in our special case of DQG (which has a bounded counit) we observe the following
   \blmma\label{defn_lemma}
   Given a (non-degenerate) coassociative $C^*$-algebra homomorphism $\al:\cla\to \clm(\cla\otimes_{\rm min}\clq)$ where $\clq$ is a $C^*$-algebraic DQG and $\cla$ is a $C^*$-algebra, the following are equivalent.
   \begin{enumerate}
       \item[(1)] $\al$ is injective.
       \item[(2)] $({\rm id}\otimes \eps)\al(a)=a$ for all $a\in \cla$.
   \end{enumerate}
   \elmma
   {\it Proof:} \\ (1) follows from (2) trivially.  Conversely, let (1) holds. Then we have an injective $C^*$-algebra homomorphism $\al:\cla\to\clm(\cla\otimes_{\rm min}\clq)$ which is coassociative, i.e. $(\al\otimes {\rm id})\al(.)=({\rm id}\otimes \Delta)\al(.)$. Define $\al_{\eps}(a):=({\rm id}\otimes \eps)\al(a)$ for all $a\in \cla$. Then, $$\al(\al_{\eps}(a))=({\rm id}\otimes \eps)(\al\otimes {\rm id})\al(a)=({\rm id}\otimes {\rm id}\otimes \eps)(\al\otimes {\rm id})\al(a)=({\rm id}\otimes {\rm id}\otimes \eps)({\rm id}\otimes \Delta)\al(a)=\al(a).$$Therefore, the injectivity of $\al$ implies $({\rm id}\otimes \eps)\al(a)=a$ for all $a\in \cla$. \qed \\ In view of the above, the condition of injectivity in the definition of a $C^*$-algebraic coaction of a $C^*$-algebraic DQG as in \cite{neshveyev_tuset}, is equivalent to the condition (2) of the preceding lemma. Moreover, as we will eventually show, the condition of span density is an automatic consequence of condition (2) of the above lemma in this case, which allows us to define the coaction of a $C^*$-algebraic DQG on a $C^*$-algebra as follows: 
   \bdfn\label{coaction_defn}
    For $C^*$-algebraic DQG  $\clq$, a coaction on a $C^*$-algebra $\cla$ is a unital $C^*$-algebra homomorphism $\delta: \cla \raro \clm(\cla \ot_{\min} \clq)$ which is coassociative, i.e. $({\rm id}\otimes \Delta)\circ\delta=(\delta\otimes {\rm id})\circ \delta$ such that $({\rm id}\otimes \epsilon)\delta(a)=a$ for all $a\in \cla$.\edfn 
   
\section{Discrete quantum automorphism group}
 The first author would like to thank Stefaan Vaes for useful discussion through private communication on the material in this subsection. This has already been introduced briefly in \cite{dg_outer}, but we explicitly describe the main result with proof in this section. We would like to mention that, although we originally intended to give a detailed proof in a preprint named `Discrete quantum group of outer automorphisms of von Neumann algebras' by the first named author, as referred in \cite{dg_outer}, we have now decided to give the proofs in this present article only, to make it self contained. We will not include these proofs in any other preprint or article written by us.
 
     Fix a unital $C^*$-algebra $\cla$ and consider a category $\widetilde{\cld_{\cla}}$ whose objects are given by $X=(\alpha, \clh)$, where $\clh$ is a finite dimensional complex Hilbert space, $\alpha : \cla \raro \cla \ot \clb(\clh)$ is a  unital $\ast$-homomorphism. Morphisms  ${\rm Mor}(X,Y)$ for  objects $X=(\alpha, \clh)$ and $Y=(\beta,  \clk)$ are given by  linear maps $T : \clh \raro \clk$  
      such that $$(1 \ot T) \alpha(a)=\beta(a) (1 \ot T) ~ \forall a \in \cla.$$   Tensor product and direct sum operations are as usual, given by: 
       $$ (\alpha, \clh) \oplus (\beta, \clk):=(\alpha \oplus \beta, \clh \oplus \clk),$$
       $$ (\alpha, \clh) \otimes (\beta, \clk):=((\alpha \ot {\rm id}) \circ \beta, \clh \ot \clk).$$ The unit object is given by $\mathds{1}:=({\rm id}_\cla, \IC)$. It is easy to see that $\mathds{1} \ot X \cong X \cong X \ot \mathds{1}$ for any object $X\in \widetilde{\cld_{\cla}}$. For any projection $P\in {\rm Mor}(X,X)$ we can define the sub-object $X_P=(\al_P,\clh_P)$ where $\al_P(a)=(\rm id_{\cla}\otimes P)\al(a)(\rm id_{\cla}\otimes \iota_{P\clh})$ where $\iota_{P\clh}:P\clh\to \clh$ is the inclusion map.  Therefore, $\widetilde{\cld_{\cla}}$ is a $C^*$ tensor category or unitary tensor category (UTC) which is closed under taking finite direct sums, tensor products and sub-objects.
    \\We call an object $(\al,\clh_{\al})\in \widetilde{\cld_{\cla}}$ to be dualizable if the following conditions hold:\\
     there are $s \in \overline{\clh} \ot \clh$, $t \in \clh \ot \overline{\clh}$ (where $\overline{\clh}$ denotes the complex conjugate Hilbert space of $\clh$), $\hat{\alpha} : \cla \raro \cla \ot \clb(\overline{\clh})$  unital $\ast$-homomorphism satisfying
     $$( (\hat{\alpha}  \ot {\rm id}) \circ \alpha(a)) (1 \ot s)=a \ot s,~~~((\alpha \ot {\rm id}) \circ \hat{\alpha}(a))(1 \ot t)=a \ot t~~~\forall a \in \cla,$$
     $$ (R^*_s \ot 1_{\overline{\clh}})(1_{\overline{\clh}} \ot R_t)=1_{\overline{\clh}},~~~~(R^*_t \ot 1_\clh)(1_\clh \ot R_s)=1_\clh,$$
      where for a vector $v $ in some finite dimensional vector space $\clk$ we denote by $R_v$ the map from $\IC$ to $\clk$ given by $z \mapsto zv $. 
\\Let us denote $\cld_\cla$ as the subcategory of all dualizable objects in $\widetilde{\cld_{\cla}}$. Therefore, invoking Theorem 2.4 of \cite{Longo_Roberts}, we have $\cld_{\cla}$ is a  rigid full subcategory of $\widetilde{\cld_{\cla}}$ which is also closed under tensor product, finite direct sum and sub-objects.    
      We have a canonical fiber functor $F$ from the UTC $\cld_\cla$ to Hilb (the UTC of finite dimensional Hilbert spaces) given by $F(\al,\clh)=\clh$ and $F(T)=T$ for any morphism $T.$ By Tannaka-Krein duality theorem this will give us a  DQG, or equivalently the dual CQG. 
      \bdfn\label{qaut_defn}We call the above discrete quantum group as the 
         universal DQG of quantum automorphism for $\cla$, denoted by $\clq_{\rm aut}( \cla)$. \edfn
      At the algebraic level, the DQG corresponds to ${\rm Nat}(F)$, the set of natural transformations from $F$ to $F$, with the algebra operation given by composition and the coalgebra structure coming from the fact that $F$ is monoidal. Indeed, any element $T \in {\rm Nat}(F)$ is uniquely determined by 
         $$\{ T_X \in \clb(\clh_X),~~ T_X T=TT_X~\forall T \in {\rm Mor}(X,X),~~~X \in {\rm Obj}(\cld_\cla) \},$$ where $X=(\alpha_X, \clh_X)$, say. This can further be identified with the subset of the above set indexed by the irreducible objects $X$ only, in which case ${\rm Mor}(X,X) \cong \IC I$, so that we have 
         $\{ T_X \in \clb(\clh_X) ,~~ X \in  {\rm Irr}(\cld_\cla) \}.$ The $C^*$-algebraic DQG $\clq_{\rm aut}(\cla)$ is given by the $C^*$-algebraic closure of $\bigoplus_{X\in {\rm Irr}(\cld_{\cla})}\clb(\clh_X)$ in $\clb(\bigoplus_{X\in {\rm Irr}(\cld_{\cla})}\clh_X$), which is nothing but the $C_0$ sequences in $\prod_{X\in {\rm Irr}(\cld_{\cla})}\clb(\clh_X)$. 
         We define $\Theta : \cla \raro \clm(\cla \ot \clq_{\rm aut}(\cla)) \cong \prod_{X \in {\rm Irr}(\cld_\cla)} \cla \ot \clb(\clh_X)$ by \begin{equation}\label{eqnforcoaction}
         \Theta(a)=\prod_{X\in \rm Irr(\cld_{\cla})}\alpha_X(a), ~~~a \in \cla.\end{equation} As each $\alpha_X$ is a unital $\ast$-homomorphism, it follows that $\Theta$ is so. Moreover, we can verify \be \label{coasso}(\Theta \ot {\rm id}) \circ \Theta=({\rm id} \ot \Delta) \circ \Theta, \ee using the definition of $\Delta$ in the Tannaka-Krein theorem through the identification of ${\rm Nat}(F \ot F)$ with ${\rm Nat}(F) \ot {\rm Nat}(F)$ and the fact that $\alpha_{X \ot Y}=(\alpha_X \ot {\rm id})\circ \alpha_Y$. 
         Let us recall the definition of the coproduct $\Delta$. The algebra $\cln:={\rm Nat}(F)$ is a multiplier Hopf algebra in the sense of Van Daele with the algebraic multiplier $\clm(\cln \ot_{\rm alg} \cln)$ being identified with a subset of the algebraic direct product of $\clb(\clh_X) \ot \clb(\clh_Y) $ with $X,Y$ varying over all objects of 
          $\cld_\cla$. Any element $W$  of this direct product can be uniquely given by its `components' indexed by $(X,Y)$, say $W_{(X,Y)} \in \clb(\clh_X) \ot \clb(\clh_Y) \cong \clb(\clh_X \ot \clh_Y)$.  In this picture, the coproduct of the multiplier Hopf algebra $\cln$ is nothing but the map which sends $T \equiv \{T_X: ~X \in {\rm Obj}(\cld_\cla) \}$ to $\Delta(T) \equiv \{ \Delta(T)_{(X,Y)}=T_{(X \ot Y)} \}.$ One can show that the restriction of this algebraic coproduct indeed maps the $C^*$-algebra $\clq_{\rm aut}(\cla)$ into $\clm(\clq_{\rm aut}(\cla) \ot_{\rm min} \clq_{\rm aut}(\cla))$.
          
         For any bounded functional $\omega $ on $\cla$, 
         write $A^\omega$ for $(\omega \ot {\rm id})(A)$ for any $A \in \clm(\cla \ot_{\rm min} \clq)$ for any other $C^*$-algebra $\clq$.
          Thus, $\Theta^{\omega}(a) \in \clq_{\rm aut}(\cla) =\prod_{X \in {\rm Irr}(\cld_\cla)} \alpha^\omega_X(a).$ From the discussion on the definition of $\Delta$ it is clear that $R:=\Delta(\Theta^\omega(a))$ is given by its `components', $R_{(X,Y)}=(\Theta^\omega(a))_{X \ot Y}=\alpha^\omega_{X \ot Y}(a)=
          (\omega \ot {\rm id} \ot {\rm id}) ((\alpha_X \ot {\rm id}) \circ  \alpha_Y(a))=((\Theta \ot {\rm id}) \circ \Theta(a))^\omega_{(X,Y)}.$ The last equality follows from the following calculation.\begin{align*}
          (\omega\otimes {\rm id}\otimes {\rm id})(((\Theta\otimes {\rm id})\circ\Theta)(a)_{(X,Y)})
          =&(\omega \otimes {\rm id}\otimes {\rm id})((\Theta\otimes {\rm id})\prod_{Z}\alpha_Z(a))_{(X,Y)} \\
          =&(\omega\otimes {\rm id}\otimes {\rm id})(\prod_{W,Z}(\alpha_W\otimes {\rm id})\circ \alpha_Z(a))_{(X,Y)}\\
          =&(\omega\otimes {\rm id}\otimes {\rm id})((\alpha_X\otimes {\rm id})\circ\alpha_Y(a)).
           \end{align*}  As $\omega$ is arbitrary, this proves (\ref{coasso}). 
    The above discussion yields the following 
    \bthm\label{Theta_coaction}
    The map $\Theta$ is a coaction of the $C^*$-algebraic DQG $\clq_{\rm aut}(\cla)$ on $\cla$.
    \ethm
    {\it Proof:} We will continue to use the notation $\Theta(a)=\prod_{X\in {\rm Irr}(\cld_{\cla})}\al_X(a)$ as in \eqref{eqnforcoaction}. The map $\Theta$ is coassociative because of the preceding discussion. Furthermore, $({\rm id}\otimes \eps)(\Theta(a))=({\rm id}\otimes \eps)(\prod_{X\in {\rm Irr}(\cld_{\cla})}\al_X(a))=\al_{\bf{1}}(a)=a.$ 
        Therefore, $\Theta$ is a $C^*$-algebraic coaction of $\clq_{\rm aut}(\cla)$ on $\cla$.
    \qed
    \bdfn\label{weak_faith}
    Let $\cls$ be a $C^*$-algebraic DQG given by the $C^*$-algebraic closure of the direct sum of full matrix algebras $\bigoplus_{i\in I}\clb(\clh_i)$ with a coaction $\delta$ on a unital $C^*$-algebra $\cla$. We call the coaction $\delta$ to be `weakly faithful' if the von Neumann algebra generated by the set $\{(\omega\otimes {\rm id})\delta(a):a\in \cla,\omega \in \cla^*\}$ inside $\clb(\oplus_{i\in I}\clh_i)$, is $\cls''.$( here $\cla^*$ denote the space of bounded linear functionals on $\cla$.
    \edfn
    
          Let us note the following:
          \bthm\label{thetaisweakfaith}
          $\Theta$ is weakly faithful in the sense that the von Neumann algebra generated by the elements of the form $\Theta^{\omega}(a)=(\omega \ot {\rm id})(\Theta(a)),$ where $a \in \cla$ and $\omega$ is a bounded linear functional on $\cls$, is the whole of $\clq_{\rm aut}''$.
          \ethm
          {\it Proof:}\\
          Consider the von Neumann subalgebra $\clw$ (say) of $\clq_{\rm aut}(\cla)^{\rm von}$ generated by $\Theta^{\omega}(\cla)$, $\omega$ bounded functional,  as described in the statement of the theorem.  Let us calculate the commutant $\clw^\prime$ in 
          $\clb(\bigoplus_X \clh_X)$. Writing a general element $T$ in this commutant in the block-operator form $(T_{(X,Y)})$ with $T_{(X,Y)} \in \clb(\clh_X, \clh_Y)$, we have the relation  $(1 \ot T_{(X,Y)})\alpha_X(a)=\alpha_Y(a) (1 \ot T_{(X,Y)})$ for all $a \in \cla$, hence $T_{(X,Y)}$ is a morphism 
           between the objects indexed by $X$ and $Y$, which are inequivalent irreducibles, so we must have $T_{(X,Y)}=0$ if $X$ is not $Y$ and otherwise a constant, say $c_X$ times $1_{\clh_X}$. This means $\clw^\prime={\rm diag}(c_X 1_{\clh_X} )_{X \in {\rm Irr}(\cld_\cla)}$.
           Taking commutant once more, we get $\clw=\clw^{\prime \prime}=\prod_X \clb(\clh_X)=\clq_{\rm aut}(\cla)''$.\qed\\
          
          Before we prove the main theorem establishing the universality of $\clq_{\rm aut}(\cla)$, we state and prove a few results. \\In that pursuit, we begin with some facts about discrete quantum groups and their unitary corepresentations. Let $\cls$ be a $C^*$-algebraic DQG  
         given by closed direct sum of finite dimensional matrix algebras $ \clb(\clh_i),~ i \in I$ ($I$ is some index set) with a coproduct $\Delta$. Let $\cls_0$ denote the dense $\ast$ subalgebra inside $\cla$ given by the algebraic direct sum of full matrix algebras, i.e. $\cls_0=\bigoplus_{i\in I}\clb(\clh_i)$. Then the dual CQG $\hat{\cls}$ has irreducible unitary mutually inequivalent corepresentations indexed by $I$. For $i\in I$, let the corresponding irreducible unitary corepresentation be $\clu_i \in \hat{\cls}\otimes \clb(\clh_i)$, which can be written in the matrix form $( {\clu_i})_{pq} $ with respect to some choice of orthonormal basis of $\clh_i$. The linear span of these `matrix coefficients' of the unitary corepresentations $\clu_i$ i.e. the span of $\{ {(\clu_i)}_{pq}, p,q=1, \ldots, d_i,~ i\in I \}$ is a dense Hopf $\ast$-algebra inside $\hat{\cls}$, denoted by $\hat{\cls}_0$. This Hopf algebra can be identified with the subspace of the dual of $\cls_0$ consisting of the `finitely supported' functionals, i.e. functionals which vanish outside the range of some $P_K$, with $K$ finite. The product is the convolution : $(\phi \cdot \psi)(a)=(\phi \ot \psi) (\Delta(a))$ for $a \in \cls_0$. Indeed, as noted in \cite{dqg}, for 
            finite subsets $K,L $ of $I$, $(P_K \ot P_L) \Delta(a) \in \cls_0$ for all $a \in \cls_0$, which makes the above product meaningful where $\phi, \psi$ are finitely supported functionals). We also have the involution $\ast$ on the space of finitely supported functionals given by $\phi^*(a)=\overline{\phi(S(a)^*)}$ for 
             $a \in \cls_0$. Let us from now on identify $\cls_0$ with the above Hopf $\ast$-algebra of (finitely supported) functionals on $\cls_0$ and also note that under this identification, ${(\clu_i)}_{pq}$ w.r.t. a choice $\{ e^p, p=1, \ldots, d_i \}$ of orthonormal basis of $\clh_i$, $i \in I$, gets identified with the functional $\phi^i_{pq}$ which is $1$ on $e^i_{pq}$ and $0$ on all other basis elements of $\cls_0$, where $e^i_{pq}(e^i_k)=\delta_{qk} e^i_p$ as before. 
             
        For any $i\in I$ we know from the general theory of CQGs, there is a unique $\overline{i}\in I$ such that $\clu_{\overline{i}}$ is equivalent to ${\clu_i}^{(c)}$, i.e. there is some positive invertible $\rho_i$ s.t. $ (1\otimes \rho_i) {\clu_i}^{(c)} (1\otimes \rho_i^{-1})$ is unitary in 
    $\hat{\cls}\otimes \clb(\overline{\clh_i})$ and there is a unitary $T_i$ from $\overline{\clh_i}$ to $\clh_{\overline{i}}$  s.t. 
    \begin{equation}
    \label{5}
    \clu_{\overline{i}}=(1\otimes  T_i\rho_i) \clu_i^{(c)} (1\otimes \rho_i^{-1}T_i^{-1}).
      \end{equation}
      
      Let $P_i$ denote the central projection corresponding to $\cls_i$, $P_K:=\sum_{i \in K} P_i$ for $K\subseteq I$. We call a set $K\subset I$ symmetric if for any $i\in K$, $\overline{i}\in K$. For any symmetric set $K$ and the unbounded closable antipode map $S$, $S(P_K)=P_K$.
         It is known (ref. e.g. \cite{dqg}) that there is a one-dimensional piece $\cls_{\mathds{1}} \cong \IC$, such that the map $\epsilon=P_{\mathds{1}}$ acts as the counit of the underlying multiplier Hopf algebra structure. In particular, 
       $(\epsilon \ot {\rm id})\Delta=({\rm id} \ot \epsilon) \Delta={\rm id}$ (here both $\Delta$ and $\epsilon$ are bounded, unital $\ast$-homomorphisms). 
     Let us now consider a unitary corepresentation $\clv$ of $\cls$ in a Hilbert space $\clk$. By definition, $\clv$ is a unitary in $\clm(\clb(\clk) \ot \cls)$  and $\clv_{(12)}\clv_{(13)}=({\rm id} \ot \Delta)(\clv)$. 
           Define $\Pi_{\clv}:\hat{\cls}\to \clb(\clk)$ as $\Pi_{\clv}(\phi)=(\rm id\otimes \phi)(\clv).$ Then we have the following:    \blmma
           \label{rep}
         $\Pi_\clv$ is a unital $\ast$-homomorphism from $\hat{\cls}_0$ to $\clb(\clk)$. Moreover, \begin{equation} \label{22} (\Pi_{\clv} \ot {\rm id})(\clu_i)=\clv^i, \end{equation} where we write $\clv \in \clm(\clb(\clk) \ot \cls)$ as $\bigoplus_{i \in I} \clv^i$ with $\clv^i \in \clb(\clk) \ot \clb(\clh_i)$.     \elmma
        {\it Proof:}
    
   (i)  First we prove that $\Pi_\clv$ is homomorphic. For $\phi, \psi$ with support say $K, L$ (finite subsets of $I$), \begin{align*} \Pi_\clv(\phi \cdot \psi)
    =& ({\rm id} \ot \phi \ot \psi) ( ( {\rm id} \ot \Delta) (\clv))\\
    =& ({\rm id} \ot \phi \ot \psi) (\clv_{(12)} \clv_{(13)}) \\
    =& ({\rm id} \ot \phi \ot \psi)(1 \ot P_K \ot P_L)(\clv_{(12)}\clv_{(13)})\\
    =& ({\rm id} \ot \phi \ot \psi)(\clv^K_{(12)} \clv^L_{(13)} )~~~({\rm where}~~\clv^K=(1 \ot P_K)(\clv),~~\clv^L=(1 \ot P_L)(\clv))\\
    =&({\rm id} \ot \phi)(\clv^K)({\rm id} \ot \psi)(\clv^L)\\
    =& \Pi_\clv(\phi)\Pi_\clv(\psi).
    \end{align*}
     Note that as $\clv^L, \clv^K$ are in the algebraic tensor product $\clb(\clk) \ot_{\rm alg} \cls$, above calculations make sense. Next, we observe that the counit map is the unit element of $\hat{\cls_0}$ and $W:=\Pi_\clv(\epsilon) $ 
    is  a unitary in $\clb(\clk)$ as $\epsilon$ is a unital $\ast$-homomorphism and $\clv$ is unitary. It also satisfies $W^2=W$ using the facts that $\clv$ is a corepresentation and $({\rm id}\ot \epsilon)\circ \Delta(\cdot)={\rm id}$ (all maps involved here are  bounded $\ast$-homomorphism). Hence $W=I$.
     
    The relation (\ref{22}) can easily be verified by evaluating $\Pi_\clv$ on $\phi^i_{pq}$ mentioned above, for any choice of orthonormal basis. We now prove that $\Pi_\clv(\phi^*)=\Pi_\clv(\phi)^*$ for any $\phi\in \hat{\cls}_0.$ 
    To this end, recall the closed unbounded antipode map $S$ and the fact noted before that $S(P_K)=P_K$ for any symmetric finite subset $K $ of $I$. As $P_K$ is a central projection, $\clv^K=(1 \ot P_K)(\clv)$ is a unitary in $\clb(\clk) \ot_{\rm alg} P_K\cls$, its inverse being $(\clv^K)^*=(\clv^*)^K$,
     and it is in the domain of $(1 \ot S)$.
    We claim that \begin{equation} \label{antipode} ({\rm id} \ot S)(\clv^K)=(\clv^*)^K.\end{equation} It suffices to verify that $({\rm id} \ot S)(\clv^K) \clv^K=1 \ot P_K$. Note that (see \cite{dqg}) for $a, b \in \cls_0$ (algebraic direct sum), $\Delta(a) (1 \ot b) \in \cls_0 \ot_{\rm alg} \cls_0$ and $m \circ (S \ot {\rm id})(\Delta(a)(1 \ot b))=\epsilon(a) b$ where $m$ is the algebraic 
       multiplication map. As $\clv^K$ and $({\rm id} \ot S)(\clv^K)$ are in $\clb(\clk) \ot_{\rm alg} \cls_0$, let us write $\clv^K= \sum_{i \in K, p,q} v^i_{pq} \ot e^i_{pq} $.  By Prop. 2.2. of Van Daele, \cite{dqg}, there is a finite subset $L$ of $I$ such that  for $i$ outside $L$ 
        $(P_K \ot P_K) (\Delta(P_i)) =0$. Writing $\clv$ as a strongly convergent sum $\clv=\sum_{i \in I} (1 \ot P_i)V $  the above implies, $(1\ot P_K \ot P_K) ({\rm id} \ot \Delta)(\clv)=(1 \ot P_K \ot P_K) ({\rm id} \ot \Delta)(\clv P_L)$. Without loss of generality, we can assume $\mathds{1} \in L$.
       Hence $\clv^K_{(12)} \clv^K_{(13)} =(1 \ot P_K \ot P_K) (\clv_{(12)} \clv_{(13)})=(1 \ot P_K \ot P_K) ({\rm id} \ot \Delta)(\clv^L)$. 
        We have $$({\rm id} \ot S)(\clv^K) \clv^K= ({\rm id}  \ot m \circ (S \ot {\rm id}) )((\clv^K)_{(12)} \clv^K_{(13)} )
       = ({\rm id } \ot m \circ (S \ot {\rm id}) )({\rm id} \ot (P_K \ot P_K) \Delta(\clv^L))$$
         Write $\clv^L=\sum_{i\in L,p,q} v^i_{pq} \ot E^i_{pq}$, we have $({\rm id } \ot m \circ (S \ot {\rm id}) )({\rm id} \ot (P_K \ot P_K) \Delta(\clv^L))=  \sum_{i\in I, p,q} v^i_{pq} \ot m \circ (S \ot {\rm id}) (P_K \ot P_K) (\Delta(E^i_{pq})).$ 
            As $S(P_K)=P_K$ is a central projection and   $m \circ (S \ot {\rm id})(\Delta(a)(1 \ot b))=\epsilon(a) b,$ for $a,b \in \cls_0$, the above expression becomes $\sum_{i\in L, p,q} v^i_{pq} \ot \epsilon(E^i_{pq}) P_K=v^e \ot P_K=1 \ot P_K.$ This proves our claimed relation (\ref{antipode}). 
             Finally, for a finitely supported functional $\phi$, assuming without loss of generality that $\phi(\cdot)=\phi(P_K \cdot)$ for some symmetric finite subset $K$, it is easy to observe from the definition of $\phi^*$ that $\phi^*(\cdot)=\phi^*(P_K \cdot)$ too, as $S(P_K)=P_K$. Hence 
               \begin{align*}
               &{\Pi_V(\phi^*)}
              = ({\rm id } \ot \phi^*)(V^K)
              = \sum_{i \in K, p,q} v^i_{pq} \overline{\phi((S(E^i_{pq})^*)}.
              \end{align*}
               But $\sum_{i\in K, p,q} v^i_{pq} \ot S(E^i_{pq}))=({\rm id} \ot S)(\clv^K)=(\clv^*)^K=(\clv^K)^*=\sum_{i \in K, p,q} (v^i_{pq})^* \ot (E^i_{pq})^*$ by (\ref{antipode}) and as $P_K$ is central. Hence 
               $\Pi_{\clv}(\phi^*)=\sum_{i \in K, p,q} (v^i_{pq})^* \overline{\phi(E^i_{pq})} =( \sum_{i \in K, p,q} v^i_{pq} \phi(E^i_{pq}))^*=(({\rm id }\ot \phi)(\clv^K))^*=\Pi_{\clv}(\phi)^*.$ 
               \qed
    
          \blmma\label{coaction_is_object}
              For any $C^*$-algebraic coaction $\delta:\cla\to \clm(\cla\otimes _{\rm min}\cls)$ of any DQG $\cls=\overline{\bigoplus_{i\in I}\clb(\clh_i)}$ on $\cla$, the  direct summand $(\delta_i,\clh_i)\in {\rm Obj}( \cld_{\cla}).$ 
          \elmma
          {\it Proof:}\\ Let $\delta:\cla\to \clm(\cla\otimes _{\rm min}\cls)$ be a $C^*$-algebraic coaction of the DQG $\cls=\overline{\bigoplus_{i\in I}\clb(\clh_i)}$ on the unital $C^*$-algebra $\cla$. Denoting the bounded counit of $\cls$ by $\epsilon$, we have from the definition of $C^*$-algebraic coaction as given in the Preliminaries, $({\rm id}\otimes \epsilon)\delta(a)=a$ for all $a\in \cla$. We can extend the $C^*$-algebraic coaction $\delta$ to a unital normal $\ast$-homomorphism $\tilde{\delta}:\tilde{\cla}\to \tilde{\cla}\overline{\otimes}\cls^{\rm von}$ where $\tilde{\cla}$ is the universal enveloping von Neumann  algebra of $\cla$ in $\clb(\clh_u)$ where $\clh_u$ is the Hilbert space corresponding to the universal representation of $\cla$ (see for reference \cite{dixmier_Cstaralgebra}) and $\cls^{\rm von}$ is the von Neumann algebraic closure of $\cls$ inside $\clb(\bigoplus_{i\in I}\clh_i)$. One can easily check that $\tilde{\delta}$ is coassociative. $\eps$ extends to a bounded unital normal $\ast$-homomorphism on $\cls^{\rm von}$ as mentioned before. Note that $({\rm id}\otimes \epsilon)\tilde{\delta}(a)=({\rm id}\otimes \eps)\delta(a)=a$ for all $a\in \cla$. By the normality of $\tilde{\delta}$ and $({\rm id}\otimes \epsilon)$ (where $\eps$ also denotes the bounded normal extension of $\eps$ on $\cls^{\rm von}$) and the density of $\cla$ in $\tilde{\cla}$ in the ultraweak topology we have $({\rm id}\otimes \eps)\tilde{\delta}(a)=a$ for all $a\in \tilde{\cla}$. Therefore, $\tilde{\delta}$ is injective too. Hence it is a von Neumann algebraic coaction in the sense of \cite{vaes_impl}. Therefore, by Proposition 3.7 and Theorem 4.4 of \cite{vaes_impl}, there exists a unitary corepresentation $\clv$ of $\cls^{\rm von}$ on $\clh_u$ such that $\tilde{\delta}(a)=Ad_{\clv}(a)$ for all $a\in \tilde{\cla}$. \\Now, from the direct sum decomposition of $\cls$, we can write $\tilde{\delta}=\prod_{i\in I}\tilde{\delta}_i$ such that $\tilde{\delta}_i=Ad_{\clv^i}$, where $\clv^i$ is the projection of $\clv$ on the $ith$ component. Therefore, $\delta_i(a)=\tilde{\delta}(a)=Ad_{\clv^i}(a)$ for all $a\in \cla$. We know from Lemma \ref{rep} there exist unitary finite dimensional corepresentations $\clu_i$ of $\hat{\cls}$ such that the map $\Pi:\hat{\cls}\to \clb(\clh_u)$ satisfying $(\Pi\otimes {\rm id})(\clu_i)=\clv^i$, is a unital $\ast$-homomorphism. Let us fix $i\in I$ and recall the operators $\rho_i$ and $T_i$ introduced in the discussion preceding Lemma \ref{rep}. Therefore, ${\clv^i}^{(c)}=({\rm id}\otimes J)({\clv^i}^*)=(\Pi\otimes {\rm id})({\rm id}\otimes J)(\clu_i^*)=(\Pi\otimes {\rm id})(\clu_i^{(c)})$, which implies,$\hat{{\clv^i}}:=(1\otimes \rho_i)\clv_i^{(c)}(1\otimes \rho_i^{(-1)})=(\Pi\otimes {\rm id})({\rm id}\otimes {\rho_i })\clu_{i}^{(c)}({\rm id}\otimes \rho_i^{(-1)})$ is unitary. Then, $\hat{\delta_i}:=Ad_{\hat{\clv^{i}}} :\cla\to \clm(\cla\otimes \cls)$ is a  unital $\ast$-homomorphism. Furthermore, taking $t=({\rm id}\otimes \rho)(\sum_{r=1}^{n}e_r\otimes\overline{e_r})$ and $s=(\rho\otimes {\rm id})(\sum_{r=1}^n\overline{e_r}\otimes e_r)$, (where $n$ is the dimension of the finite dimensional Hilbert space $\clh_i$) one can easily check that $(\delta_i,\clh_i)$ and $(\hat{\delta_i},\overline{\clh_i})$ are dual to each other in the category $\cld_{\cla}$. 
          \qed \\
         We also have the following 
         \bcrlre\label{span_density}
         For any $C^*$-algebraic coaction $\delta:\cla\to \clm(\cla\otimes _{\rm min}\cls)$, we have the span density condition, i.e. the linear span of $\delta(\cla)({\rm id}_{\cla}\otimes_{\rm min}\cls)$ is norm dense in $\cla\otimes _{\rm min}\cls$.
         \ecrlre
         {\it Proof:}\\
         Given any $C^*$-algebraic coaction $\delta:\cla\to \clm(\cla\otimes _{\rm min}\cls)$, of the DQG $\cls=\overline{\bigoplus_{i\in I}\clb(\clh_i)}$, we have from the previous lemma, $(\delta_i,\clh_i)$ is an object of $\cld_{\cla}$ with a dual object $(\hat{\delta},\overline{\clh_i})$. Let $\{e_i:i\in 1,2,\cdots,n\}$ and $\{\overline{e_p}:p\in 1,2,\cdots,n\}$ be orthonormal bases of the Hilbert spaces $\clh_i$ and $\overline{\clh_i}$ respectively where $n$ is the dimension of $\clh_i$. Recall the ket-bra notations where for any vectors $\omega,\eta\in \clh$, we define the operator $|\omega\rangle\langle\eta|\in \clb(\clh)$ by $|\omega\rangle\langle\eta|(\zeta):=\langle\eta,\zeta\rangle\omega.$ Using this notation, we have the well-known isomorphism $\Phi:\clh_i\otimes \overline{\clh_i}\ni e_p\otimes \overline{e_q}\mapsto |e_p\rangle\langle e_q|\in \clb
         (\clh_i)$. In particular, under this isomorphism the vector $\xi:=\sum_{r=1}^ne_r\otimes\overline{e_r}$ gets mapped to $1_{\clb(\clh_i)}.$ We use this fact in the following computation, where using Sweedler notation we write
        for any $a\in \cla$ and $\xi=\sum_{r=1}^ne_r\otimes \overline{e_r}$, $\delta_i(a)=\delta_i(a)_{(1)}\otimes \delta_i(a)_{(2)}\text{ and }\hat{\delta_i}(a)=\hat{\delta_i}(a)_{(1)}\otimes \hat{\delta_i}(a)_{(2)}.$ Then we have 
    \begin{align*}
        &(\delta_i\otimes {\rm id})\hat{\delta_i}(a)(1\otimes\xi)=a\otimes \xi\\
        \text{or, }&(\delta_i(\hat{\delta_i}(a)_{(1)})\otimes \hat{\delta_i}(a)_{(2)})(1\otimes \xi)=a\otimes \xi\\
        \text{or, }&\sum_{r=1}^n\delta(\hat{\delta}(a)_{(1)})(1\otimes |e_r\rangle\langle\hat{\delta_i}(a)_{(2)}(\overline{e_r})|)(1\otimes P_i)=a\otimes \sum_{r=1}^n|e_r\rangle\langle e_r|~\text{ where }P_i=1_{\clb(\clh_i)}\in \cls\\
        \text{or, }&\sum_{r=1}^n\delta(\hat{\delta}(a)_{(1)})(1\otimes |e_r\rangle\langle\hat{\delta_i}(a)_{(2)}(\overline{e_r})|)(1\otimes P_i)=a\otimes 1_{\clb(\clh_i)}~\text{ under the identification } \Phi
    \end{align*} which implies that $\delta(\cla)(1\otimes _{\rm min}\cls)$ is linearly dense in $\cla\otimes_{\rm min} \cls.$ 
         \qed

          \bthm
          \label{qaut_univ}
          (i) For any $C^*$-algebraic DQG $\cls$ with a coaction $\gamma : \cla \raro \clm(\cla \ot_{\rm min} \cls)$ there is a unique non degenerate $\ast$-homomorphism $\hat{\gamma}$ from $\clq_{\rm aut}(\cla)$ to $\cls$ such that $({\rm id} \ot \hat{\gamma})\circ \Theta=\gamma$. Moreover, $\hat{\gamma}$ is a DQG morphism.\\
          (ii) If, furthermore, $\gamma$ is weakly faithful then $\hat{\gamma}$ is surjective.\\
          (iii) $(\clq_{\rm aut}(\cla), \Theta)$ is the (unique upto isomorphism) universal object in the category with objects $(\clg, \beta)$ where $\clg$ is a $C^*$-algebraic DQG and $\beta : \cla \raro \cla \ot_{\rm min} \clg$ a $C^*$-algebraic coaction and morphisms from $(\clg_1, \beta_1)$ to $(\clg_2, \beta_2)$ being DQG morphisms $\eta : \clg_1 \raro \clg_2$ such that 
          $({\rm id}\otimes \eta)\circ\beta_1=\beta_2$.
          \ethm
          {\it Proof:}
          \begin{enumerate}
              \item [(i)] Writing $\cls$ as a closed direct sum of finite dimensional matrix algebras, say $\cls=\overline{\bigoplus_{i \in S}} \cls_i$, with $\cls_i \cong \clb(\clh_i)$ for some finite dimensional $\clh_i$, we can write $\gamma=\prod_i \gamma_i$. We know from Lemma \ref{coaction_is_object}, for each $i\in I$ $(\clh_i, \gamma_i)$ is an object in $\cld_\cla$ and hence we can decompose it into irreducibles, each of which is isomorphic with some $\alpha_X$ for $X \in {\rm Irr}(\cld_\cla)$, which is a direct summand of $\clq_{\rm aut}(\cla)$. For $i \in I$, let $m^X_i$ be the multiplicity of $\alpha_X$ in $\gamma_i$, i.e. $\gamma_i$ is unitarily equivalent to the finite direct sum of $\alpha_X$ repeated $m^X_i$ times,
         for those $X$ with $m^X_i$ nonzero. That is, we have a unitary $T^X_i :\clh_i \raro \bigoplus_X ( \clh_X \ot \IC^{m^X_i})$ (with the obvious convention that $\IC^{(0)}=(0)$) such that $\gamma_i(\cdot)=(1\ot (T^X_i )^* \left( \bigoplus_X (\alpha_X(\cdot) \ot 1_{m^X_i})  \right) (1 \ot T^X_i)$.  
          Now, define $\hat{\gamma}=\oplus_X \hat{\gamma}_X$, where $\hat{\gamma}_X : \clb(\clh_X) \raro \cls$  is given by $\hat{\gamma}_X(T)=\bigoplus_{i \in I}(T^X_i)^*(T\otimes {\rm id}) P^X_i T^X_i $, where $P^X_i$ denotes the projection onto the summand $\clh_X \ot \IC^{m^X_i} $ in $\bigoplus_{X\in {\rm Irr}(\cld_{\cla})}\clh_X\otimes \mathbb{C}^{m_i^X}.$ The fact that $\hat{\gamma}$ is nondegenerate follows easily from its definition. Recall the well-known fact that every non-degenerate representation of any $\ast$-subalgebra $\clc$ of $\clk(\clh)$ for some Hilbert space $\clh$ extends to a normal representation of $\clc''\subseteq \clb(\clh)$. Hence $\hat{\gamma}$ extends to a normal $\ast$-homomorphism from $\clq_{\rm aut}^{\rm von}(\cla)\to \cls$. Now we show that $\hat{\gamma}$ is indeed a DQG morphism between the $(\clq_{\rm aut}(\cla),\Delta)\to(\cls,\Delta_{\cls})$. Since $\gamma$ is a a coaction of $\cls$ on $\cla$, we have \begin{align*}
              &(\gamma\otimes {\rm id})\circ\gamma=({\rm id}\otimes \Delta_{\cls})\circ \gamma\\
              &\text{or, }(({\rm id}\otimes\hat{\gamma})\circ \Theta \otimes {\rm id})\circ ({\rm id}\otimes \hat{\gamma})\circ\Theta=({\rm id}\otimes \Delta_{\cls})\circ ({\rm id}\otimes \hat{\gamma})\circ\Theta\\
              &\text{or, }({\rm id}\otimes \hat{\gamma}\otimes\hat{\gamma})\circ(\Theta\otimes {\rm id})\circ\Theta=({\rm id}\otimes \Delta_{\cls})\circ({\rm id}\otimes \hat{\gamma})\circ\Theta\\
              &\text{or, }({\rm id}\otimes \hat{\gamma}\otimes\hat{\gamma})\circ ({\rm id}\otimes \Delta)\circ\Theta=({\rm id}\otimes \Delta_{\cls})\circ({\rm id}\otimes \hat{\gamma})\circ\Theta \\
              &\text{or, } (\hat{\gamma}\otimes\hat{\gamma})\circ\Delta(\omega\otimes {\rm id})\Theta(a)=\Delta_{\cls}\circ({\rm id}\otimes \hat{\gamma})(\omega\otimes {\rm id})(\Theta)(a)\text{ for all a}\in \cla\text{ and }\omega\in {\cla^*}.
              \end{align*} 
              As the set $\{(\omega\otimes {\rm id})(\Theta(a)):a\in \cla,\omega\in \cla^*\}$ is weakly dense in $\clq_{\rm aut}^{\rm von}(\cla)$, therefore, invoking Theorem \ref{thetaisweakfaith}, we can infer that $(\hat{\gamma}\otimes \hat{\gamma})\circ\Delta(a)=\Delta_{\cls}\otimes \hat{\gamma}(a)$ for all $a\in \oplus_{X\in {\rm Irr}
          (\cld_{\cla})}\clb(\clh_X)$. Hence they agree on the DQG $\clq_{\rm aut}^{\rm von}(\cla)$ and therefore on the $C^*$-algebraic DQG $\clq_{\rm aut}(\cla)$.
          \item[(ii)]  If $\gamma$ is weakly faithful, then for each $i$, $(\omega \ot {\rm id})\circ \gamma_i$ generates $\clb(\clh_i)$, so $(\gamma_i(\cla))^\prime$ must be $\IC$, hence  $(\clh_i, \gamma_i)$ is irreducible as an object of $\cld_\cla$. Then $\gamma_i=\alpha_{X_i}$ for some $X_i \in {\rm Irr}(\cld_\cla)$. 
          It is clear that $\hat{\gamma}_X$ is the projection to $\clb(\clh_i)$ for the unique  $i$ s.t. $X=X_i$ (if any) and $0$ if there is no such $i$. Hence $\hat{\gamma}$ is surjective.
       \item[(iii)] It follows from (i). Uniqueness upto isomorphism is a consequence of the general property of universal object in a category.  \end{enumerate} \qed\\
          \brmrk\label{classical_setup}
          We would like to remark that in the classical version, we get the maximal discrete quantum subgroup of $\clq_{\rm aut}(\cla)$ which is commutative as a $C^*$-algebra, i.e. isomorphic to $C_0(\Gamma)$ where $\Gamma\cong{\rm Aut}(\cla)_d$ which is ${\rm Aut}(\cla)$ with the  the discrete topology. 
          \ermrk

          We conclude this subsection with a brief discussion of a quantum subgroup of $\clq_{\rm aut}(\cla)$, namely the Kac type quantum automorphism group $\clq_{\rm aut}^{\rm kac}(\cla)$. The definition of Kac type quantum groups has already been recalled in the Preliminaries section. 
          \brmrk\label{remarkon Kac}
          Our main result is a first step in the direction of proving a discrete quantum group analogue of various non-existence results proved by the first author and his collaborators in the context of compact quantum groups. Kac type discrete or compact quantum groups form a simpler yet important subclass of quantum groups and as of now we have restricted ourselves to this subclass. Moreover, it is known that any CQG having a faithful coaction on $C(X)$ for a compact space $X$ must be of Kac type (for reference see \cite{huang_huichi}). Although we could not prove any such result for DQG coaction, we believe that our result can be proved in the full generality for any discrete quantum group.  
         \ermrk
         In order to obtain this subgroup, we consider a full subcategory $\cld^{\rm kac}_\cla$ of $\cld_\cla$ by restricting to the objects $(\alpha, \clh)$ for which one can choose $s=\sum_{i=1}^n \overline{e_i} \ot e_i,~~t=\sum_{i=1}^n e_i \ot \overline{e_i}$ for some orthonormal basis $\{ e_1, \ldots, e_n\}$ of $\clh$, with $\{ \overline{e_i} \}$ being the corresponding orthonormal basis in $\overline{\clh}$. Restricting the fiber functor to the subcategory $\cld^{\rm kac}_\cla$, we get $\clq_{\rm aut}^{\rm kac}(\cla).$
          \\Moreover, as in the case of $\clq_{\rm aut}(\cla)$ we can define $\Theta^{\rm kac}:\cla \to \cla\otimes \clq_{\rm aut}^{\rm kac}(\cla)$ and prove the analogue of Theorem \ref{qaut_univ} stated below, following the same line of arguments of the proof of \ref{qaut_univ}.
          \bthm
          \label{qaut_univ_kac}
          (i) $\Theta^{\rm kac} $ is a weakly faithful coaction on $ \cla $ by $ \clq_{\rm aut}^{\rm kac}(\cla)$ .\\
          (ii) For any Kac type $C^*$-algebraic DQG $\cls$ with a coaction $\gamma : \cla \raro \clm( \cla \ot_{\rm min} \cls)$ there is a unique nondegenerate $\ast$-homomorphism $\hat{\gamma}$ from $\clq_{\rm aut}^{\rm kac}(\cla)$ to $\cls$ such that $({\rm id} \ot \hat{\gamma})\circ \Theta=\gamma$.\\
          (iii) If, furthermore, $\gamma$ is weakly faithful then $\hat{\gamma}$ is surjective.\\
          \ethm
          {\it Proof:}
          \begin{enumerate}
              \item[(i)] The proof follows by a verbatim adaptation of the proof of Theorem \ref{thetaisweakfaith}.
              \item[(ii)] Let $(\cls,\Delta_{\cls})$ be a Kac type discrete quantum group with a coaction $\gamma:\cla\to \clm(\cla\otimes\cls)$. Recall from the proof of Lemma \ref{coaction_is_object}, we can write $\gamma_i(a)=Ad_{\clv^i}(a)$ for all $a\in \cla$. Further, for each $i\in I$, there exists corepresentation $\clu_i$ of $\hat{\cla}$ such that  
              $$
                  {\clv^i}^{(c)}:=({\rm id}\otimes J)({\clv^i}^*)
                  =({\rm id}\otimes J)(\Pi\otimes {\rm id})(\clu_i^*)
                  =(\Pi\otimes {\rm id})({\rm id}\otimes J)(\clu_i^*)
                  =(\Pi\otimes {\rm id})(\clu_i^{(c)}). $$
            Now, the uni-modularity of $\cls$ implies the unitarity of $\clu_i^{(c)}$ which further implies ${\clv^i}^{(c)}$ is unitary. 

              Therefore, we can easily see that $(Ad_{{\clv^i}^{(c)}},\overline{\clh_i})\in \cld^{\rm kac}_\cla$ forms a conjugate object to $\gamma_i$ under the choice of $t=\sum_{i=1}^ne_i\otimes \overline{e_i}$ and $s=\sum_{i=1}^n\overline{e_i}\otimes e_i.$.
             \item[(iii)] The proof again follows easily from the proof of part (ii) of Theorem \ref{qaut_univ}.
          \end{enumerate}
          \qed

\section{DQG Coaction on $C(S^1)$}
In this main section, our aim is to prove that there is no genuine Kac type $C^*$-algebraic DQG with a weakly faithful, linear coaction on $C(S^1)$. Let us denote by $Z$ the canonical generator of $C(S^1)$, i.e. $Z(t)=t$ for all $t \in S^1$. $\overline{Z}$ denotes the complex conjugate of $Z$. 

\bdfn
A co-action $\beta : C(S^1) \raro \clm( C(S^1) \ot \clq)$ for some DQG $\clq$ is called linear if $\beta(1)=1 \ot 1_{\clm(\clq)}$ and the linear span $\cll$ (say) of $\{ Z,\overline{Z} \} $ is invariant under the action of $\beta$, in the sense that it is mapped into $\cll \ot \clq$. 
\edfn
\brmrk\label{remarkonlinear}
We would like to mention here that the reason behind considering linear coactions is that any isometric coaction of a compact quantum group on $C(S^1)$ becomes linear as discussed in \cite{bhowmick_goswami}. It is also straightforward to see that if we do put a similar definition of isometric coactions of DQGs, the arguments of \cite{bhowmick_goswami} will go through and we can conclude that isometric coactions of DQGs are linear too. Hence Theorems \ref{main_theorem} and \ref{main_corollary} naturally hold for an isometric coactions as well. 
\ermrk
\bthm 
\label{lin_univ}
Let $\cla=C(S^1)$. Let $\cld^{\rm kac, lin}_\cla$ be the full subcategory of $\cld^{\rm kac}_\cla$ consisting of those objects $(\alpha, \clh)$ for which there is a choice of $\hat{\alpha}$ such that both $\alpha$ and $\hat{\alpha}$ maps the span $\cll$  of $\{Z, \overline{Z} \}$ into $\cll \ot \clb(\clh)$ and $\cll \ot \clb(\overline{\clh} )$ respectively. Then  we have the following:\\
(i) $\cld^{\rm kac, lin} _\cla$ is a UTC and the restriction of $F$ to this subcategory gives a DQG $\clq_{\rm aut}^{\rm kac, lin}(\cla)$.\\
(ii) There is a weakly faithful, linear coaction, say $\Theta^{\rm kac, lin} $ by $\clq_{\rm aut}^{\rm kac, lin}(\cla)$ on $\cla$.\\
(iii) An analogue of statement (ii) and (iii) of Theorem \ref{qaut_univ} hold true  by replacing $\clq_{\rm aut}(\cla)$ by $\clq_{\rm aut}^{\rm kac, lin}(\cla)$.

\ethm
{\it Proof:}\\
The proof can be done by a verbatim adaptation of the proof of Theorem \ref{qaut_univ} and is hence omitted.\qed

 The following theorem plays a crucial role in the derivation of our main result.
\bthm\label{main_theorem}
 $\clq_{\rm aut}^{\rm kac, lin} (C(S^1)) $ is commutative as a $C^*$-algebra, i.e. isomorphic with $C_0(\Gamma)$ for the discrete group $\Gamma$ which is the isometry group of the circle with discrete topology, i.e. $\Gamma\cong S^1_d\rtimes \mathbb{Z}_2$ where $S^1_d$ is the group $S^1$ equipped with the discrete topology and the semidirect product is with respect to the canonical $\mathbb{Z}_2$ action where the non trivial element of $\mathbb{Z}_2$ maps $z\mapsto \overline{z}$. 
 \ethm
 {\it Proof:}\\
 As is evident from the discussion subsequent to Definition \ref{qaut_defn} we can write $\Theta^{\rm kac, lin} $ as a direct product of object $\alpha_i : C(S^1) \raro C(S^1) \ot M_{n_i} $, for some positive integer $n_i$, where $\alpha_i, i \in I$ (say) is an enumeration of (mutually inequivalent) 
  irreducible objects of $\cld^{\rm kac, lin}_\cla$. Thus, it suffices to prove that for any object $(\alpha, \clh)$ of this category, the image of $C(S^1)$ under $\alpha$ is contained in $C(S^1) \ot \clv$ for some {\it commutative} subalgebra $\clv$ of $\clb(\clh)$. Let us prove this by some explicit detailed computations.
 
Before we begin let us note that $\{1,Z,\overline{Z}, Z^2, \overline{Z}^2\}$, where $Z:S^1\to \mathbb{C}$ is the map sending $z\mapsto z$ and $\overline{Z}:S^1\to \mathbb{C}$ is the map sending $z\mapsto \overline{z}$ is a linearly independent subset of $C(S^1)$. Let $(\alpha, \clh= \IC^n)$ be an object of $\cld^{\rm kac, lin}_\cla$ with the corresponding $\hat{\alpha}$. Let us fix an orthonormal basis $\{ e_i \}$ of $\clh = \IC^{n}$ and the corresponding basis $\{ \overline{e_i} \}$ of $\overline{\clh}$ and view both $\clb(\clh)$ and $\clb(\overline{\clh} )$ as $M_n(\IC)$ 
 with respect to these bases. With this, and also using linearity, let us write $$\mathcal{\alpha} \colon C(S^1) \to C(S^1)\otimes M_n~~\text{ and } ~~\mathcal{\hat{\alpha}} \colon C(S^1) \to C(S^1) \otimes M_n$$ satisfying
\begin{align}\label{eqn1}
\mathcal{\alpha}(Z)= Z\otimes A + \overline{Z}\otimes B \\ 
\label{eqn2}
\text{and }\mathcal{\hat{\alpha}}(Z)= Z\otimes C + \overline{Z}\otimes D.
\end{align}
Then we have 
\begin{align}
\label{three}
((\mathcal{\hat{\alpha}}\otimes {\rm id})\circ\mathcal{\alpha}(Z))(1\otimes  \sum_{i=1}^n (\overline{e_{i}}\otimes e_{i}))= Z\otimes \sum_{i=1}^n\overline{e_i}\otimes e_{i},\\
\label{four}
((\mathcal{\hat{\alpha}}\otimes {\rm id})\circ\mathcal{\alpha}(\overline{Z}))(1\otimes  \sum_{i=1}^n (\overline{e_{i}}\otimes e_{i}))= \overline{Z}\otimes \sum_{i=1}^n\overline{e_i}\otimes e_{i},\\
\label{five}
((\mathcal{\alpha}\otimes {\rm id})\circ\mathcal{\hat{\alpha}}(Z))(1\otimes \sum_{i=1}^n(e_{i}\otimes \overline{e_{i}}))= Z\otimes \sum_{i=1}^n e_{i}\otimes \overline{e_{i}}\\
\label{six}
\text{and }((\mathcal{\alpha}\otimes {\rm id})\circ\mathcal{\hat{\alpha}}(\overline{Z}))(1\otimes \sum_{i=1}^n(e_{i}\otimes \overline{e_{i}}))= \overline{Z}\otimes \sum_{i=1}^n e_{i}\otimes \overline{e_{i}}
\end{align} where $1$ is the unit element of $C(S^1)$ i.e., $1(t)=1$ for all $t\in S^1$.
Since $\alpha$ is a unital $\ast$-homomorphism, hence we have
$
\alpha(1)= 1\otimes I_{n\times n}$, which implies, 
\begin{align*}
I_{n\times n}&=\al(Z\overline{Z})
=\alpha(Z){\alpha(Z)}^*
\end{align*}
 Expanding $\alpha(Z)$ and using the linear independence of $\{1,Z,\overline{Z},Z^2,\overline{Z}^2\}$ we get
\begin{align}
\label{seven}
&AA^* + BB^*= I_{n\times n}~~\text{ and }~~
AB^*=0_{n\times n}=BA^*.
\end{align}
Similarly, since $\hat{\alpha}$ is also a unital $\ast$-homomorphism hence we have 
\begin{align}
\label{eqn10}
&CC^* + DD^*=I_{n\times n}~~\text{ and }~~CD^*= 0_{n\times n}=DC^*.
\end{align}
From equation (\ref{three}) we have 
\begin{align*}
Z\otimes \sum_{i=1}^n\overline{e_i}\otimes e_{i}
&=((\mathcal{\hat{\alpha}}\otimes {\rm id})\circ\mathcal{\alpha}(Z))(1\otimes  \sum_{i=1}^n (\overline{e_{i}}\otimes e_{i}))\\
&=(\hat{\alpha}\otimes {\rm id})(Z\otimes A + \overline{Z}\otimes B)(1\otimes \sum_{i=1}^n\overline{e_{i}}\otimes e_{i})\\
&=(\hat{\alpha}(Z)\otimes A + \hat{\alpha}(\overline{Z})\otimes B)(1\otimes \sum_{i=1}^n\overline{e_{i}}\otimes e_{i})\\
&=(Z\otimes C\otimes A + \overline{Z}\otimes D \otimes A + \overline{Z}\otimes C^*\otimes B + Z\otimes D^* \otimes B)(1\otimes \sum_{i=1}^n\overline{e_{i}}\otimes e_{i})\\
&=Z\otimes \sum_{i=1}^nC(\overline{e_i})\otimes A(e_{i}) + \overline{Z}\otimes \sum_{i=1}^nD(\overline{e_i})\otimes A(e_i) + \overline{Z}\otimes \sum_{i=1}^nC^*(\overline{e_{i}})\otimes B(e_{i}) \\
&+ Z\otimes \sum_{i=1}^nD^*(\overline{e_{i}})\otimes B(e_{i}),
\end{align*}
which further implies,
$$Z[\sum_{i=1}^nC(\overline{e_{i}})\otimes A(e_{i}) + \sum_{i=1}^nD^*(\overline{e_{i}})\otimes B(e_{i}) - \sum_{i=1}^n\overline{e_{i}}\otimes e_{i}]+ \overline{Z}[\sum_{i=1}^nD(\overline{e_{i}})\otimes A(e_{i}) + \sum_{i=1}^nC^*(\overline{e_{i}})\otimes B(e_{i})]$$
$= 0_{n\times n}.$\\
Hence, we have 
\begin{align}
\label{eqn13}
&\sum_{i=1}^nC(\overline{e_{i}})\otimes A(e_{i}) + \sum_{i=1}^nD^*(\overline{e_{i}})\otimes B(e_{i}) - \sum_{i=1}^n\overline{e_{i}}\otimes e_{i}= 0_{n\times n}\\
\label{eqn14}
\text{and }&\sum_{i=1}^nD(\overline{e_{i}})\otimes A(e_{i}) + \sum_{i=1}^nC^*(\overline{e_{i}})\otimes B(e_{i})= 0_{n\times n}.
\end{align}
Implementing the same procedure as above on equation \eqref{four}, we have 
\begin{align}\label{eqn15}
&\sum_{i=1}^nC^*(\overline{e_i})\otimes A^*(e_{i}) + \sum_{i=1}^nD(\overline{e_{i}})\otimes B^*(e_{i}) - \sum_{i=1}^n\overline{e_{i}}\otimes e_{i} = 0_{n\times n}
\\ 
\label{eqn16}\text{ and } &\sum_{i=1}^nD^*(\overline{e_{i}})\otimes A^*(e_{i}) + \sum_{i=1}^nC(\overline{e_i})\otimes B^*(e_{i})= 0_{n\times n}.
\end{align}
Similarly, from equation \eqref{five} we have 
\begin{align}\label{eqn17}
&\sum_{i=1}^nA(e_i)\otimes C(\overline{e_i}) + \sum_{i=1}^nB^*(e_i)\otimes D(\overline{e_i})- \sum_{i=1}^ne_i\otimes \overline{e_i}= 0_{n\times n}\\
\label{eqn18} \text{ and }&\sum_{i=1}^nB(e_i)\otimes C(\overline{e_i}) + \sum_{i=}^nA^*(e_i)\otimes D(\overline{e_i})= 0_{n\times n}.
\end{align}
Finally, equation \eqref{six} translates into the following equations, namely
\begin{align}\label{eqn19}
&\sum_{i=1}^nA^*(e_i)\otimes C^*(\overline{e_i}) + \sum_{i=1}^nB(e_i)\otimes D^*(\overline{e_i}) - \sum_{i=1}^ne_i\otimes \overline{e_i}= 0_{n\times n}\\
\label{eqn20}\text{ and }&\sum_{i=1}^nB^*(e_i)\otimes C^*(\overline{e_i}) + \sum_{i=1}^nA(e_i)\otimes D^*(\overline{e_i})=0_{n\times n}.
\end{align}
Let us write the matrices A, B, C and D in $ M_n(\IC)$ as follows:$$
A=(a_{ij})_{i,j};~~
B=(b_{ij})_{i,j};~~
C=(c_{ij})_{i,j};~~
D=(d_{ij})_{i,j}
$$
where $i,j \in \{1,2,\ldots n\}.$ 
We now introduce $X,Y\in \mathrm{M}_n(\mathbb{C})$ such that $C\overline{e_i}=\overline{Xe_i}$ and $D\overline{e_i}=\overline{Ye_i}$. Then, under the identification of $\clh\otimes\overline{\clh}\cong\clb(\clh)\cong \overline{\clh}\otimes \clh$, we obtain from \begin{align*}
    \sum_{i=1}^nC(\overline{e_i})\otimes A(e_i)
    &=\sum_{i=1}^n\overline{Xe_i}\otimes A(e_i)
    =\sum_{i=1}^n\sum_{j=1}^nX_{ji}^*\overline{e_j}\otimes \sum_{k=1}^nA_{ki}e_k\\
    &=\sum_{i,j,k=1}^nX_{ji}^*A_{ki}\overline{e_j}\otimes e_k
    =\sum_{i,j,k=1}^nA_{ki}X^*_{ij}\overline{e_j}\otimes e_k
    =\sum_{j,k=1}^nAX^*_{kj}\overline{e_j}\otimes e_k\\
    &=AX^*.
\end{align*}  Similarly, we have $$\sum_{i=1}^nD^*(\overline{e_i})\otimes B(e_i)=\sum_{j,k=1}^nBY_{kj}\overline{e_j}\otimes e_k=BY. $$
Therefore equation \eqref{eqn13} under the above identifications yields
\begin{align}\label{new1}
   AX^*+BY=1 
\end{align}
\noindent Similarly, from equation \eqref{eqn14} we have, 
\begin{align*}
AY^*+BX=0
\end{align*}
Following the above procedure, we obtain the following identities from equations \eqref{eqn13}-\eqref{eqn20}:

\begin{align*}
    &AX^*+BY=1,~~~~~~AY^*+BX=0,\\
    &X^*A+YB=1,~~~~~~Y^*A+XB=0,\\
    &AX^*+B^*Y^*=1,~~~~BX^*+A^*Y^*=0,\\
    &A^*X+BY=1,~~~~~B^*X+AY=0.
\end{align*}
Since  $\overline{Z}Z=1$ (the constant 1 function on $S^1$), hence
going along the same way as for equation \eqref{seven} we have 
\begin{align}
   \label{eqn49} A^*A + B^*B = I_{n\times n}\text { and }B^*A = 0=A^*B\\
    \label{eqn52}C^*C + D^*D =I_{n\times n}\text{ and }D^*C = 0=C^*D .
\end{align}
Multiply equation \eqref{seven} on the right by $B$ to get 
   $~AA^*B + BB^*B = B 
    \text{ which gives }BB^*B = B $
i.e.  B is a partial isometry.
\\Multiply \eqref{eqn49} on the right by $A^*$ to get
$A^*AA^* + B^*BA^* = A^*
\text{ which gives }A^*AA^* = A^*$ i.e. 
 A is a partial isometry.
Similarly $X$ and $Y$ are also partial isometries.
\\Hence there exist projections $P_1,P_2,Q_1$ and $Q_2$ such that $A^*A =P_1, B^*B=1-P_1, AA^*=P_2, BB^*=1-P_2, X^*X=Q_1, Y^*Y=1-Q_1, XX^*=Q_2\text{ and }YY^*=1-Q_2.$ 
\\Let us define $U:= A+B$ and $V:=X+Y$ which are both unitary. Observe that $A-B$ is also a unitary.\\
Further we have, from equations \eqref{seven} and \eqref{eqn49},
$$AB^*=0=BA^*\text{ and }B^*A=0=A^*B.$$ Therefore, multiplying $AY^*+BX=0$ by $BB^*$ from the left we obtain $BX=0$. Therefore, $AY^*=0$. Similarly, we obtain from the above relations $$Y^*A=XB=BX^*=A^*Y^*=B^*X=AY=0.$$
Therefore, we have $$(A+B)(X+Y^*)^*=AX^*+BY+AY+BX^*=AX^*+BY=1.$$ Sine, $A+B$ is unitary, hence we have $A+B=X+Y^*.$ Going along the same lines, we have $A-B=X-Y^*.$ Hence, $A=X$ and $B=Y^*$. Thus, $AY^*=0=Y^*A$ implies $AB=0=BA.$ 
    Finally, we claim that $A\text{ and }B$ are normal. This will complete the proof that the $C^*$-algebra generated by $A\text{ and }B$ is a commutative one. 
    The claim holds due to the following argument.
    \\As $A=UP_1$, $B=U(1-P_1)$, it suffices to prove $U$ and $P_1$ commute. As $U$ is unitary, $AB=0$ implies $$ UP_1U(1-P_1)=0 \text{ or, } P_1U(1-P_1)U^*=0 \text{ or, } P_1=P_1UP_1U^* \text{ or, } P_1U=P_1UP_1.$$ On the other hand, $BA=0 \text{ or, } (1-P_1) UP_1U^*=0 \text{ or, } UP_1=P_1UP_1$. Thus, $UP_1=P_1U$. \\This implies, $AA^*=UP_1P_1^*U^*=UP_1U^*=P_1UU^*=P_1$ and $A^*A=P_1U^*UP_1=P_1$ thus proving the normality of $A$ and therefore the normality of $B$ as well. 
    Thus, we can conclude that $\clq_{\rm aut}^{\rm kac, lin}(C(S^1))''$ is a commutative von Neumann algebra which implies the commutativity of $\clq_{\rm aut}^{\rm kac, lin}(C(S^1))$ as a $C^*$-algebra as well. Therefore, there exists a discrete group $\Gamma$ such that $\clq_{\rm aut}^{\rm kac,lin}(C(S^1))\cong C_0(\Gamma)$. In fact, this discrete group is nothing but the universal object in the category of groups acting linearly on $S^1$. Let $G$ be any group acting linearly on $S^1$. Then for any $g\in G$ and $z\in S^1$, we have the map $(g,z)\mapsto a_gz+b_g\overline{z}$. Now, $1=\|a_gz+b_g\overline{z}\|^2=\|a_g\|^2+a_g\overline{b_g}+\overline{a_g}b_g+\|b_g\|^2$ such that $a_gb_g=0$ ( by our observation in the above computations). Therefore, the only choice of $(a_g,b_g)$ is $\{(a_g,b_g):\|a_g\|=1\text { and }b_g=0\}$ or $\{(a_g,b_g):a_g=0\text{ and }\|b_g\|=1\}$. Therefore for any group $G$ acting linearly on $S^1$, we have the injective map $\Psi:G\to S^1\rtimes \mathbb{Z}_2$ given by $g\mapsto (a_g,\overline{0})$ if $b_g=0$ 
    and $g\mapsto (b_g,\overline{1})$ if $a_g=0$. Furthermore, one can easily have a linear action of $S^1\rtimes \mathbb{Z}_2$ on $S^1$. Therefore, we infer that $\Gamma$ is isomorphic to $S^1\rtimes \mathbb{Z}_2$ endowed with the discrete topology, i.e. $\Gamma\cong S^1_d\rtimes \mathbb{Z}_2.$
    We end with the main result of our paper which is a direct corollary of the previous theorem.
    \brmrk\label{diff_in_DQGcase}
    Let us explain why the ideas and techniques used in the earlier work on quantum isometry groups in the framework of CQG do not apply to the present context. For a CQG coaction, one can average over the Haar state for a CQG coaction, which is not possible for a DQG coaction as the Haar weight is not a bounded everywhere defined functional. To be more precise, given a faithful linear coaction of a CQG on $C(S^1)$, we can obtain a faithful state on $C(S^1)$  which is invariant under the CQG coaction (by averaging technique), so that the restriction of the coaction on span $\{Z,\overline{Z}\}$ becomes a two-dimensional unitary corepresentation. Then the arguments of Theorem 2.2 of \cite{bhowmick_goswami} would apply to prove that the CQG is commutative i.e. not genuine. However, we cannot use such arguments in the framework of DQG. \ermrk
    \bthm\label{main_corollary}
    Let $(\cls,\gamma)$ be a Kac type $C^*$-algebraic DQG coacting on $C(S^1)$ such that the coaction $\gamma$ is weakly faithful and linear in the sense described previously. Then $\cls$ is also commutative as a $C^*$-algebra and hence is classical, i.e. isomorphic to $C_0(\Gamma)$ for some discrete group $\Gamma$.
    \ethm
    {\it Proof:}\\ Let $(\cls,\gamma)$ be a Kac type DQG such that the coaction of $\cls$ on $C(S^1)$ is weakly faithful and linear. We know from Theorem \ref{lin_univ} that there exists an unique surjective nondegenerate $\ast$-homomorphism $\hat{\gamma}$ from $\clq_{\rm aut}^{\rm kac, lin}$ to $\cls$. Thus the surjectivity of the map $\hat{\gamma}$ and commutatitivity of $\clq_{\rm aut}^{\rm kac, lin}(\cla)$ implies that $\cls$ is also commutative as a $C^*$-algebra and hence is classical.\qed
    \\We end this paper with the following 
    \brmrk
    It remains an interesting open question whether like the non-existence of genuine weakly faithful linear $C^*$-algebraic DQG coaction on $C(S^1)$, we can prove non-existence of smooth DQG coaction for $C(M)$, for a more general compact connected smooth manifold $M$, as an exact analogue of \cite{goswami_joardar},\cite{goswami}. 
    
    \ermrk
    \textbf{ Acknowledgement} The first author would like to thank Stefaan Vaes for his suggestion of the definition of $\clq_{\rm aut}(\cla)$ through private communication. The authors would like to express their deep gratitude to the  referees whose suggestions led to significant modification of this paper and notable shortening of the proof of Theorem \ref{main_theorem}.




\bibliographystyle{amsalpha} 
\bibliography{references} 
\end{document}